\documentclass[12pt]{article}
\usepackage[utf8]{inputenc}
\usepackage{amsmath}
\usepackage{amsfonts}
\usepackage{amssymb}
\usepackage{amsthm}
\usepackage{relsize}
\usepackage{array}
\usepackage{multirow}
\usepackage{tabu}
\usepackage{soul}
\usepackage{graphicx}
\usepackage{epigraph}
\usepackage{float}
\usepackage[english]{babel}
\usepackage[usenames, dvipsnames]{color}
\usepackage{fancyhdr}
\usepackage[Sonny]{fncychap}
\usepackage{tikz-cd}
\usepackage{tkz-euclide}
\usepackage{sseq}
\usepackage{newfloat}
\usepackage{amsmath,mathtools}
\usepackage{pgfplots}
\usepackage{geometry}
 \usepackage{sectsty}
\sectionfont{\fontsize{14}{15}\selectfont}

\usepackage[all]{xy}
\DeclareFloatingEnvironment[fileext=dia]{diagram}

\theoremstyle{definition}

\usepackage[most]{tcolorbox}

\makeatletter
\newtcolorbox{note}[1][]{%
	breakable,
	enhanced jigsaw, 
	borderline west={3pt}{0pt}{black!10!white}, 
	borderline south={1pt}{0pt}{black!10!white}, 
	borderline east={1pt}{0pt}{black!10!white},
	borderline north={1pt}{0pt}{black!10!white},
	sharp corners, 
	boxrule=0pt, 
	attach title to upper, 
	left=0pt,
	right=0pt,
	top=0pt,
	bottom=0pt,
	boxsep=5pt,
	colback=white,
	frame hidden,
	#1
}
\makeatother

\makeatletter
\newtcolorbox{note1}[1][]{%
	breakable,
	enhanced jigsaw, 
	sharp corners, 
	boxrule=0pt, 
	attach title to upper, 
	fontupper=\linespread{1.1}\fontfamily{qpl}\selectfont,
	fontlower=\linespread{1.1}\fontfamily{qpl}\selectfont, 
	left=0pt,
	right=0pt,
	top=0pt,
	bottom=0pt,
	boxsep=3pt,
	colback=green!3!white,
	frame hidden,
	before skip=10pt plus 2pt,after skip=10pt plus 2pt,
	#1
}
\makeatother

\makeatletter
\newcommand\tabfill[1]{%
	\dimen@\linewidth
	\advance\dimen@\@totalleftmargin
	\advance\dimen@-\dimen\@curtab
	\parbox[t]\dimen@{#1\ifhmode\strut\fi}%
}
\makeatother

\usepackage{tabularx}
\usepackage{imakeidx}
\usepackage{hyperref}
\hypersetup{colorlinks={true},linkcolor={blue},citecolor={blue},urlcolor={blue}, linktoc=all}  
 
 \usepackage[capitalize, nameinlink]{cleveref}
 \crefdefaultlabelformat{#2\textbf{#1}#3}
 \crefname{figure}{Figure}{Figures} 
 \Crefname{figure}{Figure}{Figures}
 \crefname{table}{Table}{Tables}
 \Crefname{table}{Table}{Tables}
 \crefname{section}{\S\hspace{-1mm}}{\S\hspace{-1mm}}
 \Crefname{section}{\S\hspace{-1mm}}{\S\hspace{-1mm}}
 \crefname{equation}{}{}
 \Crefname{equation}{}{}
 \crefname{example}{Geometric Pattern}{Geometric Patterns} 
 \Crefname{example}{Geometric Pattern}{Geometric Patterns}

 \usepackage{footnote}

\begin{document}

\title{\textbf{Bisection of Trapezoids in Elamite Mathematics}}

\author{Nasser Heydari and Kazuo Muroi}

\maketitle

\begin{abstract}
 The bisection of trapezoids  by transversal lines has many examples in Babylonian mathematics. In this article, we study a similar problem in   Elamite mathematics,   inscribed on a clay tablet held in the collection of the Louvre Museum and thought to date from between 1894--1595 BC.   We   seek to demonstrate  that this problem is different from   typical Babylonian problems about bisecting trapezoids by transversal lines. We also identify some of  the possible mathematical ideas underlying  this problem and the    innovative approach that might have motivated its design.\\
 
  \noindent
 \textbf{Acknowledgment.} The authors wish to express their deep gratitude and respect to the memory of late Professor T. Taniguchi. He was the first  to realize the   creativity behind the   problem examined in this article  and  made mention of    it to the second author.   
\end{abstract}

\section{Introduction}
This tablet is one of 26 excavated from Susa in  southwest Iran by French archaeologists in 1933. The texts of all the Susa mathematical tablets (henceforth \textbf{SMT}) along with their interpretations were first published in 1961 (see \cite{BR61}). In Elamite mathematics,   two problems relating to the transversal bisectors of trapezoids   are found in \textbf{SMT No.\,23} and four problems in \textbf{SMT No.\,26}.

Bisection of trapezoids is often categorized  among the inheritance problems  found in Babylonian mathematics (see \cite{Fri05-1,Fri07-1,Fri07-2,FA16,Hyp02,Mur01-1,Oss18}, for examples of such problems). In these problems, a trapezoidal land area  is  divided equally between two brothers. One of the many ways to achieve this task is to use a line parallel to the bases of the trapezoid (a transversal). This specific method of dividing   trapezoids has been the source of great interest in Babylonian mathematical texts and the  many problems  regarding this issue   found there. Recently, it has   been suggested the Babylonians could compute the  total distance of Jupiter's travel  along the ecliptic during a certain interval of time  from the area of a trapezoidal figure representing the planet's changing daily displacement along the ecliptic. Moreover, the time when Jupiter reaches half the total distance  may be computed by bisecting the trapezoid into two smaller ones of equal areas by a transversal bisector  (see \cite{Oss16,Oss18}). It should be noted that the  reasons why Babylonians favored  transversal bisectors remain  unknown to us, but one may suggest both mathematical and nonmathematical rationales.

 However, the fourth   problem of   \textbf{SMT No.\,26} deals with this issue   differently. The problem considers  the use of transversal strips--not lines--to partition a trapezoid and tries to determine the transversal strip (if any) that  bisects  the trapezoid. Treating   such a problem in this way involves finding  the natural solutions to a quadratic equation whose solvability depends entirely on the values of the bases and the number of the strips. As discussed later, this problem  rarely has a solution and   finding a correct one  requires checking   many different values. It is a surprise that the Susa scribe of  \textbf{SMT No.\,26} seems to have successfully solved the problem.

This partly broken tablet\footnote{The reader can see the new  photos of this tablet on the website of the Louvre's collection. Please see \url{https://collections.louvre.fr/ark:/53355/cl010186436} for obverse  and   reverse.}   contains four problems each of which deals with the partition of a right trapezoid. The terminology of this text differs somewhat from other Susa mathematical texts  in  both its vocabulary and concise expression.

\section{Bisection of Geometric Figures}
\begin{figure}[H]
	\centering
	\includegraphics[scale=1]{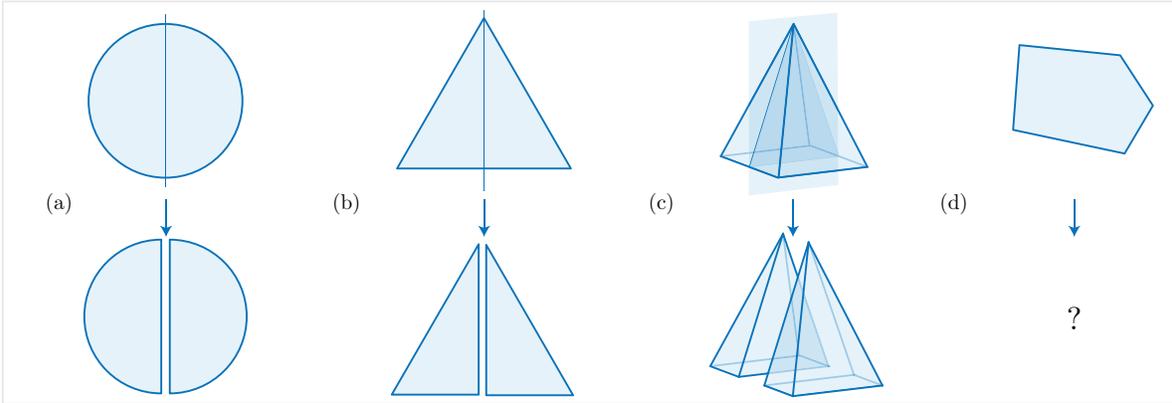}
	\caption{Bisection of different shapes}
	\label{Figure1}
\end{figure}

Bisection of a geometric shape\footnote{We only consider shapes which have bounded areas or volumes and are made up of one piece such as polygons and polyhedra. In   topology terms, these shapes are called \textit{compact} and \textit{connected} subsets of the plane and the 3-space. See \cite{Rud76} for definitions of compactness and connectedness.}  usually involves the division of the shape into two  congruent parts. For two dimensional  figures this is   done by a line   while for three dimensional shapes   a plane is used. The bisecting  lines  or  planes are    called   \textit{bisectors}. It should be noted that under the congruency condition the bisection of some figures might not be possible.  \cref{Figure1} shows the bisection of a few geometric figures with their congruent parts and bisectors.  The pentagon in case (d)  can not be bisected into congruent parts at all. 

\begin{figure}[H]
	\centering
	\includegraphics[scale=1]{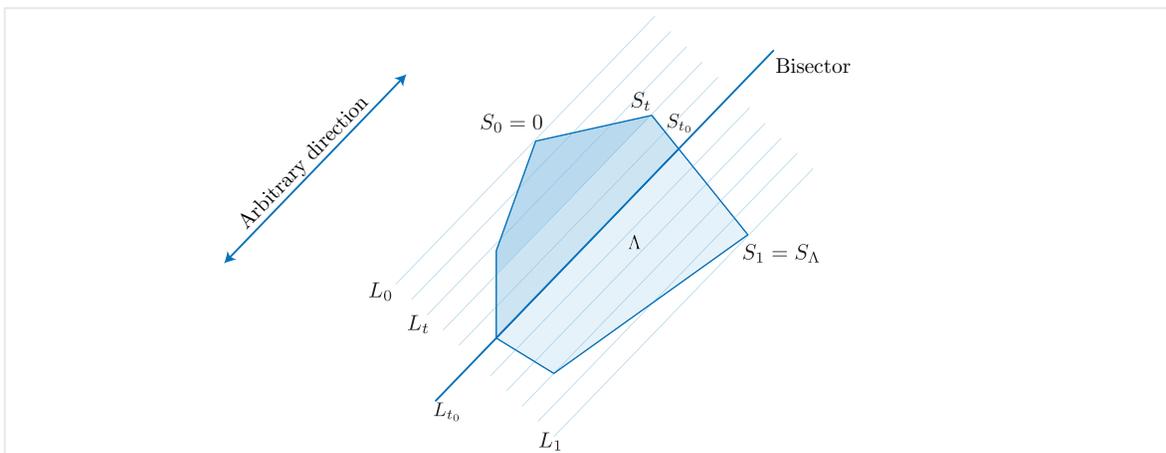}
	\caption{Bisection of a plane figure in an arbitrary direction}
	\label{Figure2}
\end{figure}

Now, let us replace the congruency condition in the definition of the bisection with a new one which only   considers  the equality of areas or volumes of the obtained parts. This concept of bisection is more general than the previous one and from now on whenever we mention the bisection of a figure, we mean the process of dividing it by a line or a plane into two parts with equal areas or volumes.   It is interesting that under this new condition any figure  can be bisected in an  infinite    number  of ways. For example, as is shown in \cref{Figure2}, for any given direction in the plane and any   plane figure  $\Lambda$, there is a line segment passing through it  and dividing it into two equal parts.

One   explanation for the above-mentioned claim can be provided by the  well-known    \textit{Intermediate  Value Theorem}\footnote{The theorem says that if $f(x)$ is a continuous function such that $f(x_1)<y<f(x_2)$, then there exists a number $x_0$ between $x_1$ and $x_2$  for which $f(x_0)=y$. For a  proof of this theorem and its applications, see \cite{Rud76}.} from mathematical analysis.      The key idea here is to consider parallel lines in the fixed direction and define a   function  on real numbers $0\leq t\leq 1$ which assigns  to each number $t$     the area $S_t$ of the part of the figure $\Lambda$  surrounded  by two parallel lines $L_0$ and $L_t$ (see \cref{Figure2}). Clearly, $S_0=0$ and $S_1=S_{\Lambda}$. The   \textit{Intermediate  Value Theorem} implies that there is a $0<t_0<1$ for which 
\[S_{t_0}=\frac{1}{2} S_{\Lambda},\]
meaning that the corresponding line $L_{t_0}$ is a bisector for the figure $\Lambda$. We should note that although this method   guarantees the existence of the bisector,   determining the equation or the length of the bisecting line segment might not be   so easy.

\begin{figure}[H]
	\centering
	\includegraphics[scale=1]{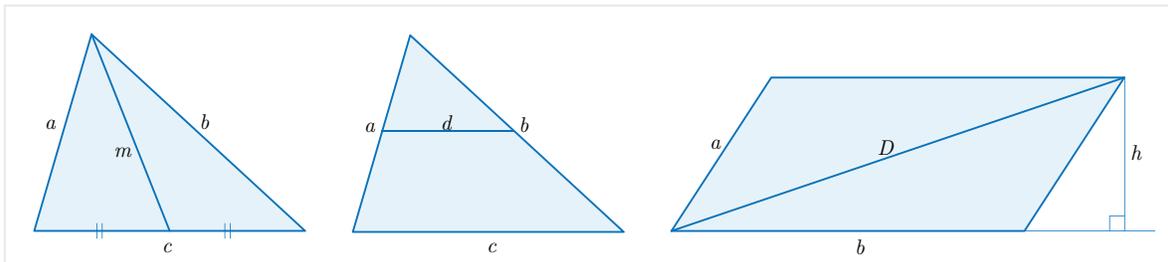}
	\caption{Bisectors of    triangles and   parallelograms}
	\label{Figure3}
\end{figure}

 One  problem in classic geometry is to compute the lengths of     bisectors of   polygons whose sides are known. For example, consider  a triangle with sides $a,b,c$.  Among all its bisectors, we can consider medians of its sides  or the transversals parallel to its  sides. In these two cases,  we can easily compute the  lengths of bisectors (see \cref{Figure3}). The length $m$ of the median     bisecting the side $c$ of the triangle   is  obtained by
\[ m =\frac{\sqrt{2a^2+2b^2-c^2}}{2}, \]
while a transversal line $d$ parallel to the side $c$  of the triangle has the length 
\[ d =\frac{\sqrt{2}c}{2}. \]

On the other hand, for a parallelogram with sides $a,b$ and height $h$, any line passing through its centroid bisects the figure. For example, two of these bisectors are the diagonals of the parallelogram whose lengths can be computed with respect to $a,b,h$. If $D$ is a diagonal of the parallelogram as in \cref{Figure3}, then
\[ D =\sqrt{a^2+b^2+2b\sqrt{a^2-h^2}}. \]

\section{Bisection of Trapezoids in Babylonian Mathematics}
In comparison with other basic geometric figures, the trapezoid has been paid  a great deal of  attention  in Babylonian mathematical tablets. In those texts, Babylonian scribes addressed many problems dealing  with bisecting a trapezoid and the length of the corresponding bisector. Of infinitely many possible bisectors   for a trapezoid, the Babylonian scribes seemed to have  considered  only  the transversal line    parallel to the two bases of the trapezoid (see \cite{Fri05-1,Fri07-1,Fri07-2,FA16,Hyp02,Mur01-1,Oss18}, for examples of such problems). 

\begin{figure}[H]
	\centering
	\includegraphics[scale=1]{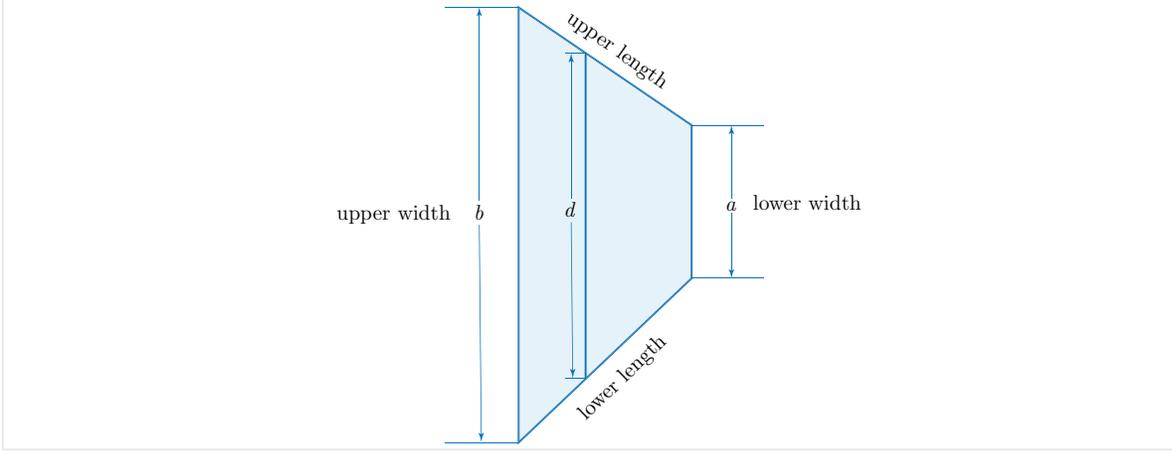}
	\caption{Transversal bisector of a trapezoid}
	\label{Figure4}
\end{figure}

 One definite advantage of choosing the transversal bisector is that its length depends only on the length of two bases  and we do not need to know the   length  of other sides or the height. In fact, as we   see shortly,  if the lengths of two bases  of a trapezoid are $a$ and $b$ (see \cref{Figure4}), then the length of the transversal bisector  $d$ is given by 
\begin{equation}\label{equ-a}
	d=\sqrt{\dfrac{a^2+b^2}{2}}.
\end{equation}

\noindent
\textcolor{Red}{\textbf{Convention.}} In Babylonian mathematical tablets, the trapezoids are usually drawn in a way that the bases are perpendicular to the horizontal direction. Because of this, the bases are called the upper and the lower widths and the other two sides are called the upper and the lower lengths (see \cref{Figure4}). 

\begin{figure}[H]
	\centering
	\includegraphics[scale=1]{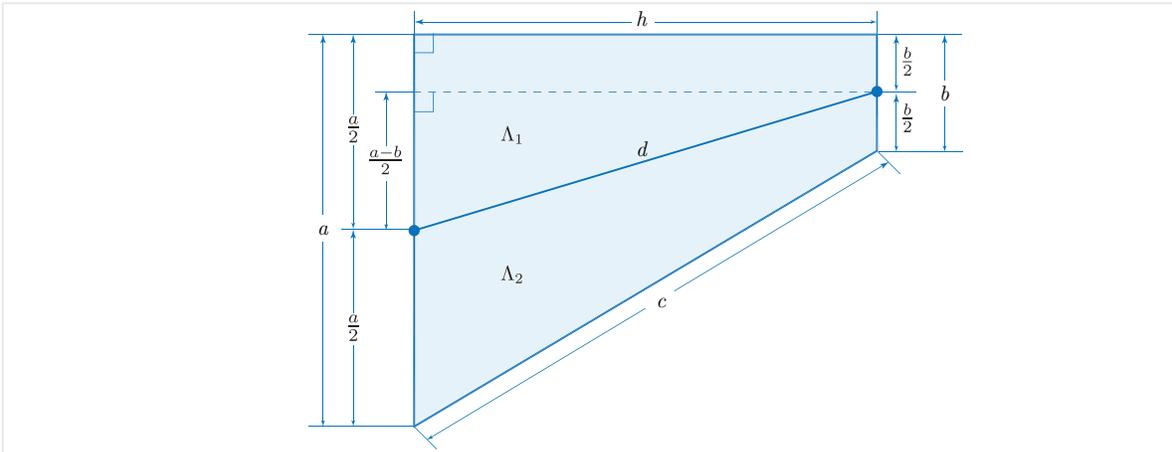}
	\caption{Bisection of a  trapezoid using midpoints of bases}
	\label{Figure5}
\end{figure}

It is  clear that formula \cref{equ-a} is independent from height of the trapezoid. This  is not true in general  for other bisectors. For example, as is evident from \cref{Figure5}, the line segment connecting the midpoints of bases divides the trapezoid into two smaller   trapezoids\index{trapezoid} $ \Lambda_1$  and $ \Lambda_2$ with equal areas, because their bases and heights have the same lengths\index{length}:
$ S_{\Lambda_1}=S_{\Lambda_2}=\frac{h}{4}(a+b)$. 
Note that the length\index{length} of the dividing    line  $d$  in \cref{Figure5} can be computed by using the Pythagorean theorem  as
\begin{equation}\label{equ-b}
	d=\sqrt{h^2+\left(\frac{a-b}{2}\right)^2} =\sqrt{c^2-3\left(\frac{a-b}{2}\right)^2}.  
\end{equation}
In both formulas of \cref{equ-b}, the value of $d$ does not depend only on the length  of bases.

\begin{figure}[H]
	\centering
	\includegraphics[scale=1]{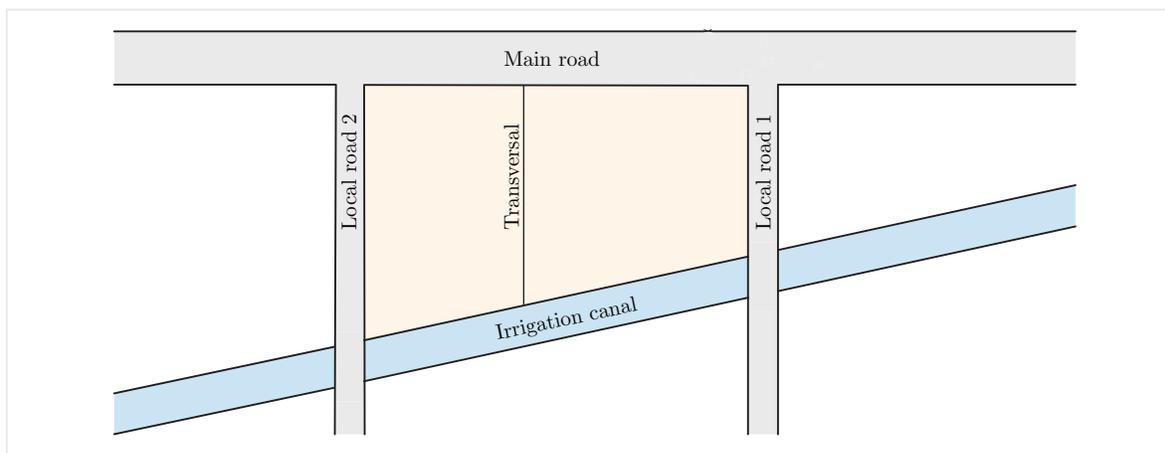}
	\caption{A trapezoidal area of farmland}
	\label{Figure6}
\end{figure}

Another reason  behind this   Babylonian choice might have been   that  this division produces two trapezoids of the same kind (for example, right trapezoids\footnote{A trapezoid having at least two right angles is called a right trapezoid.} are divided into right trapezoids).  
Aside from mathematical advantages,   social justice might also be reflected in this  method  which  utilizes  mathematical skills   to solve   real-life problems. Historically, division of an inheritance amongst members of a family has always been a   delicate issue. An unsatisfactory outcome can start a bitter family feud.  Consider a trapezoidal area of farmland as shown in \cref{Figure6} and assume that  a land surveyor is asked to divide it between two brothers equally. As is shown in the figure, three sides of the land are surrounded by roads and the fourth one is bounded by an irrigation canal.

In addition to dividing the land into two parts with equal areas, the land surveyor ought to include two important factors in his calculations: (1) direct access to the main road and  (2) direct access to the irrigation canal.  These two factors, even nowadays, are  the most decisive ones in increasing the value of any farmland.  The transversal bisecting method is the one giving  these two  advantages to both divided parts. As we can see from \cref{Figure6},   the    transversal line bisects the farmland into two equal parts each of which can directly access both  the irrigation canal and the main road.  Moreover,  the shape of the two parts are the same (both right trapezoids), a characteristic that   other division methods (such a  midpoint method)  may not possess. By dividing the land in this way, the surveyor is able to kill  two birds with one stone! He  succeeds in  his task  both   mathematically  and also  from a social justice point of view. It should be noted that   in other methods, the owner of one part of the divided land  might lose direct access to the main road or the irrigation canal.

\subsubsection*{Proof of the Babylonian formula}
\begin{figure}[H]
	\centering
	\includegraphics[scale=1]{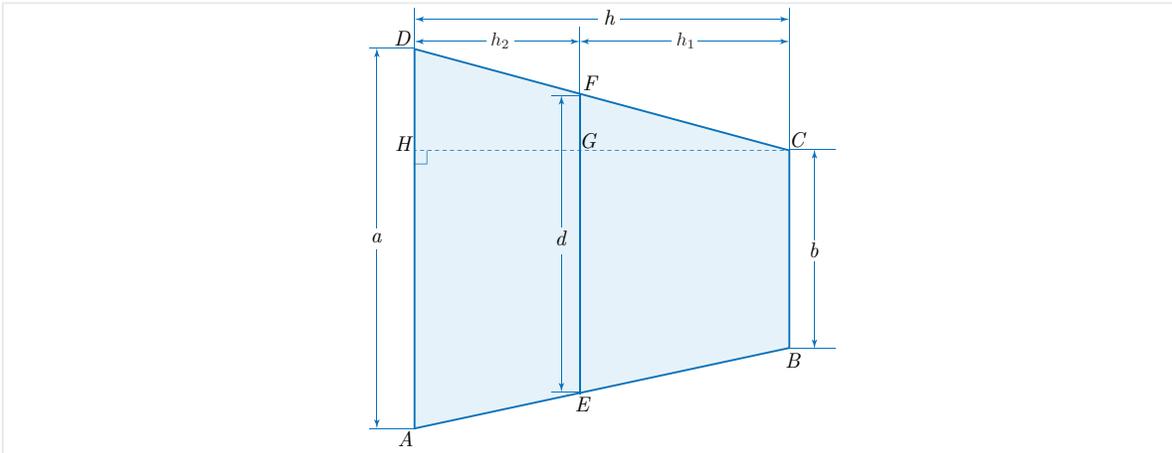}
	\caption{Bisecting transversal of a trapezoid}
	\label{Figure7}
\end{figure}

Now, we prove   the  Babylonian formula\index{Babylonian formula for the bisecting transversal of a trapezoid} \cref{equ-a}  for the bisecting transversal of   a  trapezoid\index{trapezoid}\index{bisection of a right trapezoid}.    Consider  a  trapezoid\index{trapezoid}   $ABCD$  the lengths of whose upper and lower  widths (bases) are  $\overline{AD}=a$ and $\overline{BC}=b$ respectively (see \cref{Figure7}). Suppose  that this trapezoid\index{trapezoid}   is divided into two equal parts by the transversal\index{transversal line} line segment $EF$ of length\index{length} $d$ which is parallel to the bases $ AD$ and $BC$.   Consider the perpendicular\index{perpendicular line} lines   from   $C$ to $AD$, the dotted lines  in \cref{Figure7},   and set
\begin{equation*} 
	\begin{cases}
		\overline{CG} =h_1,\\
		\overline{HG}= h_2,\\
		\overline{CH}= h.
	\end{cases}
\end{equation*}
Clearly, $h=h_1+h_2$ and 
\begin{equation*} 
	\begin{dcases}
		S_{AEFD}=\frac{(a+d)h_2}{2},\\
		S_{EBCF}=\frac{(b+d)h_1}{2}.
	\end{dcases}
\end{equation*}
Since $S_{AEFD}= \frac{1}{2}S_{ABCD}$, we have 
\begin{equation*} 
	\dfrac{(a+d)h_2}{2}=\dfrac{(a+b)h}{4} 
\end{equation*}
which implies that
\begin{equation}\label{equ-c}
	\dfrac{h_2}{h}=\dfrac{a+b}{2(a+d)}.
\end{equation}
Similarly, it follows from $S_{EBCF}= \frac{1}{2}S_{ABCD}$ that
\begin{equation}\label{equ-d}
	\dfrac{h_1}{h}=\dfrac{a+b}{2(b+d)}.
\end{equation}

Now, by using \cref{equ-c} and \cref{equ-d}  we can write
\begin{align*}
	&~~   \dfrac{h_1+h_2}{h}=1 \\
	\Longrightarrow~~&~~  \dfrac{h_1}{h}+\dfrac{h_2}{h}=1 \\  
	\Longrightarrow~~&~~  \dfrac{a+b}{2(b+d)}+\dfrac{a+b}{2(a+d)}=1 \\  
	\Longrightarrow~~&~~   \frac{a+b}{2}\left(\frac{1}{b+d}+\frac{1}{a+d}\right)=1 \\ 
	\Longrightarrow~~&~~    \frac{a+b+2d}{(a+d)(b+d)} =\frac{2}{a+b} \\ 
	\Longrightarrow~~&~~    a^2+b^2+2ab+2d(a+b)  =2d^2+2(a+b)d+2ab \\ 
	\Longrightarrow~~&~~   2d^2= a^2+b^2  \\ 
	\Longrightarrow~~&~~   d =\sqrt{\dfrac{a^2+b^2}{2}}              
\end{align*}
which completes the proof of \cref{equ-a}.

\section{Bisection of Trapezoids in the Susa Mathematical Texts}
Of the six problems in the  \textbf{SMT}  regarding the bisection of trapezoids  five of them might be termed standard problems  whose calculations involve the transversal bisectors of trapezoids and the application  of the Babylonian formula \cref{equ-a}. However, the fourth problem in  \textbf{SMT No.\,26} is different from these typical problems and instead of a line, the Susa scribe has considered a transversal party wall to bisect a trapezoid. We now consider this fourth problem  of  \textbf{SMT No.\,26}  and seek to illuminate its  mathematical significance.  

\subsection{Bisection of Trapezoids by Transversal Strips}
\begin{figure}[H]
	\centering
	\includegraphics[scale=1]{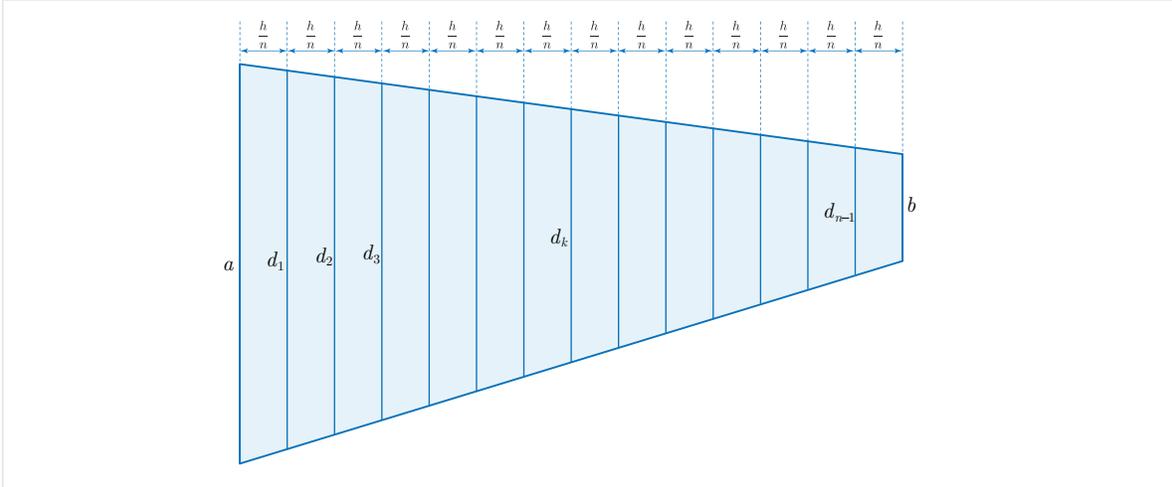}
	\caption{Partition of a trapezoid by transversal lines}
	\label{Figure8}
\end{figure}

  Consider a  trapezoid of bases $a >b$ and height  $h$. For any natural number  $n>2$, we can divide its  height  $h$ into $n$ equal parts of lengths $\frac{h}{n}$ by drawing $n-1$ vertical line segments  of lengths  $d_1,d_2,\ldots,d_{n-1}$, which are parallel to the two bases. Note that  these transversal lines are labeled from left to right.  In this case, we have   partitioned   our main   trapezoid into  $n$ smaller    trapezoids   $\Lambda_1,\Lambda_2,\ldots,\Lambda_n $. If we set $d_0 = a$ and $d_n = b$, then  the  bases   of $\Lambda_k $ are  of lengths  $  d_{k-1}$ and $ d_k  $, for all $k=1,2,\ldots,n$.  Note also that    the  heights of all trapezoids $\Lambda_k $  have  the common value $\frac{h}{n}$ (see \cref{Figure8}).

\begin{figure}[H]
 		\centering
 		\includegraphics[scale=1]{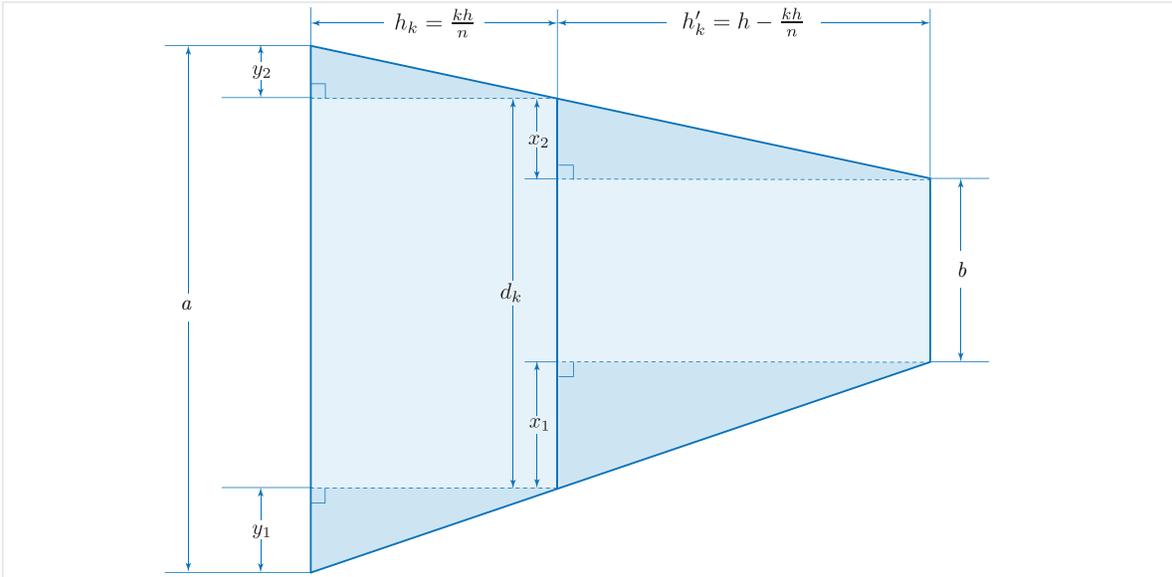}
 		\caption{Similar triangles in a  trapezoid divided by a   transversal line}
 		\label{Figure9}
 	\end{figure}

It is easy to compute the length  $d_k$  of each transversal line (bases of the small right trapezoid   $\Lambda_k$) by using the similarity  of right triangles. This situation is shown in \cref{Figure9} in which $d_k$ is the common base of two trapezoids with heights $h_k=\frac{kh}{n}$ and $h'_k=h-\frac{kh}{n}$.

It is clear from  \cref{Figure9} that $d_k=x_1+x_2+b$ and $a=y_1+y_2+d_k$.  	Since the transversal $d_k$ is parallel to bases, the two heavy-shaded right triangles in the lower part of \cref{Figure9} are similar. So, 
\[ \frac{x_1}{h'_k}=\frac{y_1}{h_k} \Longrightarrow  \frac{y_1}{\frac{kh}{n}}=\frac{x_1}{h-\frac{kh}{n}}  \]
which easily implies that
\begin{equation}\label{equ-e}
	y_1=\left(\frac{k}{n-k}\right)x_1.
\end{equation} 
 	The same reasoning for the two right triangles in the upper part of  \cref{Figure9}  also implies that
 	\begin{equation}\label{equ-f}
 		y_2=\left(\frac{k}{n-k}\right)x_2.
 	\end{equation} 
 Next, we can use \cref{equ-e} and  \cref{equ-f} and  write as follow:
 \begin{align*}
 		&~~   d_k+y_1+y_2=a \\
 		\Longrightarrow~~&~~  d_k+ \left(\frac{k}{n-k}\right)x_1 + \left(\frac{k}{n-k}\right)x_2=a \\  
 		\Longrightarrow~~&~~   d_k+ \left(\frac{k}{n-k}\right)(x_1 + x_2)=a  \\  
 		\Longrightarrow~~&~~    d_k+ \left(\frac{k}{n-k}\right)(d_k-b)=a  \\  
 		\Longrightarrow~~&~~      \left(\frac{n}{n-k}\right)d_k =a + \left(\frac{k}{n-k}\right)b \\  
 		\Longrightarrow~~&~~    d_k  = \left(\frac{n-k}{n}\right)a+  \left(\frac{k}{n}\right)b.            
 \end{align*}
Thus, we get
 \begin{equation}\label{equ-g}
	d_k=\left(1-\frac{k}{n}\right)a+\left(\frac{k}{n}\right)b,\ \ \ \text{for}\ \ \ 1\leq k \leq n-1.
\end{equation}  	
 	
\begin{figure}[H]
	\centering
	\includegraphics[scale=1]{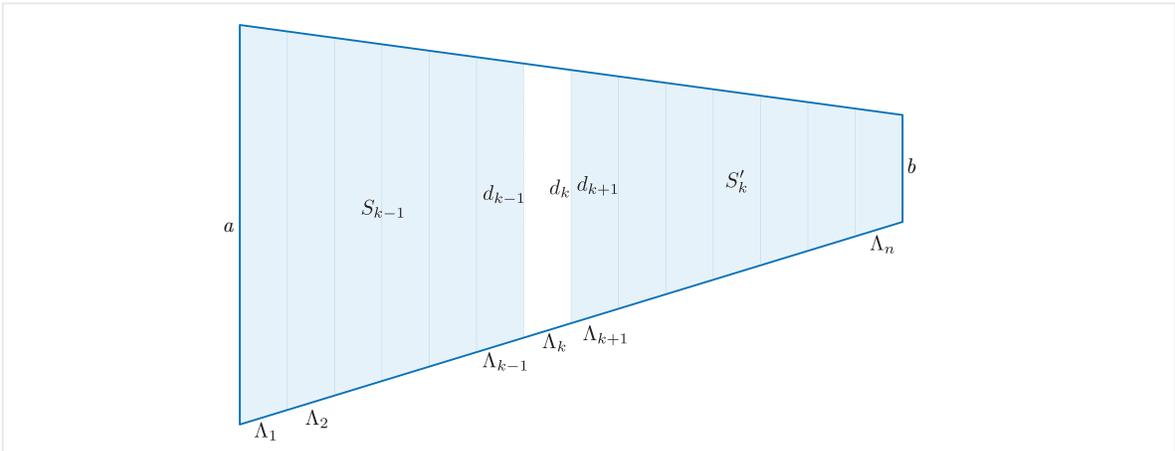}
	\caption{Grouping transversal strips}
	\label{Figure10}
\end{figure}

We raise the question that under what conditions would one of the  trapezoids $\Lambda_k$ play  the role of a transversal bisector? In another words, is it possible to choose a  	$\Lambda_k$ such that the total areas of trapezoids on its left-hand side and those on its right-hand side are equal? This situation is shown in \cref{Figure10} in which $S_{k-1}$ is the total areas of trapezoids $\Lambda_1,\Lambda_2,\ldots,\Lambda_{k-1} $ and $S'_k$ is that of trapezoids   $\Lambda_{k+1},\ldots, \Lambda_n$.

To answer this question, one needs to compute the values of $S_k$ and $S'_k$ and find the index $k$ such that    $ S_{k-1}=S'_k$.  It is clear that
\[ S_k=\frac{(a+d_1)h}{2n}+\frac{(d_1+d_2)h}{2n}+\cdots+\frac{(d_{k-1}+d_k)h}{2n}.  \]
A simple calculation  using \cref{equ-g} implies  that
\begin{equation}\label{equ-h}
S_k= \frac{kh}{2n}  \Bigg(\left(2-\frac{k}{n}\right)a+\left(\frac{k}{n}\right)b\Bigg),\ \ \ \text{for}\ \ \ 1\leq k \leq n-1,
\end{equation}
which is exactly the area  of a   trapezoid  with bases $a,d_k$ and height $h_k=\frac{kh}{n}$ (see \cref{Figure10}). Also  note that the total area of   $n-k$ small   trapezoids   $\Lambda_{k+1},\Lambda_{k+2},\ldots,\Lambda_n  $   is 
\begin{equation}\label{equ-i}
	S'_k= \frac{h\left(a+b\right)}{2} -S_k,\ \ \ \text{for}\ \ \ 1\leq k \leq n-1,
\end{equation} 
which is the area of a   trapezoid  with bases $d_k,b$ and height\index{height} $h'_k=h-\frac{kh}{n} $.

Assume that one of these $n$  small   trapezoids, say $\Lambda_{k_0} $,  is  the bisecting party wall (note that $1<k_0<n$). So,    the remaining  $n-1$  small   trapezoids   belong to  two groups:\\
(1) the ones on the left-hand side  of the party wall\index{party wall}, i.e.,   $\Lambda_{1},\Lambda_{2},\ldots,\Lambda_{k_0-1}   $,   and\\
(2)  the ones on the right-hand side of the party wall\index{party wall}, i.e, $\Lambda_{k_0+1},\Lambda_{k_0+2},\ldots,\Lambda_{n}  $. \\
Since we require  the total   areas  of these two groups to be equal,       the problem boils down to finding the possible value of $k_0$ for which   $S_{k_0-1}=S'_{k_0}$.  It  immediately follows from \cref{equ-i}   that 
\[  S_{k_0}+S_{k_0-1} = \frac{h\left(a+b\right)}{2}.  \]
By using  \cref{equ-h}  in the last equality and doing some calculations, we obtain the following quadratic equation  with respect to the unknown variable $k_0$:
\begin{equation}\label{equ-SMT26-pa}
	2(a-b)k_0^2-(4na-2b+2a)k_0+n^2(a+b)+2na+a-b=0.
\end{equation}

Clearly, the solvability of the last quadratic equation deeply depends on the values of $a,b$ and $n$. In fact, the  discriminant $\Delta$ of this equation is 
\begin{equation}\label{equ-j}
	\Delta= 4(2n^2-1)a^2+4(2n^2-1)b^2+8ab 
\end{equation}
which is always positive, because  $n>1$. Since    only  natural values between 1 and $n-1$  for $k_0$ are allowed, we definitely require the discriminant to be a perfect square which obviously is not true in general.   Therefore, one may say that  while this problem   does not  always have a solution,  under what conditions might  one   find a correct solution?

Let    $a=br$, where $r>1$. In that case, the value of discriminant $\Delta$ and $k_0$  are  
\begin{equation}\label{equ-k}
	\Delta= 4b^2\left((2n^2-1)(r^2+1)+2r\right)
\end{equation}
and
\begin{equation}\label{equ-l}
	k_0= \dfrac{(2n+1)r-1\pm \sqrt{(2n^2-1)(r^2+1)+2r}}{2(r-1)}.
\end{equation}

For any fixed $r>1$, we can try different values for $n$ and see what happen. If one check values  $r=2,3,4,\ldots,20$ and $n=3,4,5,\ldots,1000$, they will obtain the following answers:
\begin{table}[H]
\begin{center}
	\begin{tabular}[h]{c|c|c|c}
$a=rb$	&$r$ & $n$& $ k_0$\\
	\hline
$a=2b$	&	2 & 37    & 16 \\
$a=3b$	&	3  & 17  & 7 \\
$a=3b$	&	3  & 305    & 117 \\
$a=4b$	&	4  & 65    & 24\\
$a=5b$	&	5  & 10    &  4  \\
$a=6b$	&	6 & 25   & 9  \\
$a=8b$	&	8  & 35   & 12  \\
$a=9b$	&	9 & 20  & 7  \\
$a=12b$	&	12	&11	&4\\
$a=13b$	&	13&	246	&78\\
$a=15b$	&	15&	511	&160\\
$a=17b$	&	17&	8&	3\\
$a=17b$	&	17&	505	&157\\
$a=18b$	&	18&	89&	28
\end{tabular}
\caption{Acceptable values for $r$ and $n$}
\label{Table1}
\end{center}
\end{table}

This table shows that among almost 20000 cases, only 14   provide  us with  correct answers. Even if $r=2,3,\ldots,211$ and $n=3,4,\ldots,1000$ which produce  209580 cases,  there are only 32 such answers! In other words, this problem rarely has a solution and its solvability  entirely depends on the values of $\frac{a}{b}$ and $n$.

\section{Fourth Problem of \textbf{SMT No.\,26}}
Unlike the typical problems on trapezoids which deal  with   transversal bisectors,   the fourth problem of \textbf{SMT No.\,26} considers   a  situation similar to the one we discussed in the previous section. The Susa scribe seems to divide an area of  trapezoidal land   between two brothers equally by using   a party wall. This problem is treated in lines 1-17 on the reverse of this tablet. We   give  both the transliteration and the translation of this text and then explain its mathematical calculations in detail.

\subsection*{Transliteration}
\begin{note1} 
	\underline{Obverse,  Lines 1-15}\\
	(L1)\hspace{2mm} \textit{\v{s}\`{a}} $\cdots $[$\cdots $ $\cdots $] $\cdots $ [$\cdots $]
	\\
	(L2)\hspace{2mm} 2,10 sag a[n-na] ugu 30 s[ag ki-ta 1,40 dirig]\\
	(L3)\hspace{2mm} igi-3,45 u\v{s} du$_8$-\textit{ma} 16 [\textit{a-na} 1,40 íl]\\
	(L4)\hspace{2mm} 26,40 a-r\'{a} 2  53,20 \textit{a-n}[\textit{a} 4,30 íl 4]\\ 
	(L5)\hspace{2mm} \textit{ta-ta-a-ar} 2,10 sag an-n[a k\'{u}-k\'{u}-\textit{ma} 4,41,40]\\
	(L6)\hspace{2mm} 4 \textit{i-na} 4,41,40 kud-\textit{ma} [41,40]\\
	\noindent\textcolor{Green}{\rule{\textwidth}{1pt}}
	(L7)\hspace{2mm} \textit{\v{s}umma} (BAD) 50 dal murub$_4$ 30 nindan sag an(sic)-na(sic) 40(?) [$\cdots $]\\
	\noindent\textcolor{Green}{\rule{\textwidth}{1pt}}
	(L8)\hspace{2mm} 50 dal murub$_4$ ugu 30 sag ki-$<$ta$>$ \textit{mi-nam} [dirig 20 dirig]\\
	(L9)\hspace{2mm} 1,20 a-r\'{a} 20  26,40 a-r\'{a} 2  5[3,20]\\
	(L10)\hspace{0mm} \textit{ta-ta-a-ar} 50 dal murub$_4$ k\'{u}-[k\'{u}-\textit{ma} 41,40]\\
	(L11)\hspace{0mm} \'{i}b-si$_8$-bi \textit{ka-bi-is} [1,4(?)]\\
	\noindent\textcolor{Green}{\rule{\textwidth}{1pt}}
	(L12)\hspace{0mm} sag-ki-gud 35 sa[g an-na] 5 [s]ag k[i-ta]
	\\
	(L13)\hspace{0mm} 2 \v{s}e\v{s}-me\v{s} \textit{mi-it-ha-ri-i\v{s} l}[\textit{i-zu-zu}]
	\\
	\noindent\textcolor{Green}{\rule{\textwidth}{1pt}}
	(L14)\hspace{0mm} 35 a-r\'{a} 35  20,25 [5 a-r\'{a} 5  25 \textit{a-na}]
	\\
	(L15)\hspace{0mm} [20,25 \textit{\c{s}}]\textit{\'{i}-ib-ma} [20,25 bar] 10,25 [\'{i}b-si$_8$-bi 25]\\
	
	\underline{Reverse,  Lines 1-17}\\
	(L1)\hspace{2mm} [$\cdots $ $\cdots $ $\cdots $ $\cdots $] 3[6(?) $\cdots $]
	\\
	(L2)\hspace{2mm} 1,40[sag an-na 20 sag] ki-ta 1 u\v{s} [$\cdots $ $\cdots $]\\
	(L3)\hspace{2mm} \textit{i-na} \textit{l}[\textit{i-bi} u]\v{s} 1 k\`{u}\v{s} 6 \v{s}u-si \'{e}-[gar$_8$ dal-ba-na \textit{e-pu-u\v{s}}]\\
	(L4)\hspace{2mm} \'{i}b-tag$_4$ \textit{a-na} 2 \v{s}e\v{s}-me\v{s} [\textit{i-di-in}]
	\\
	\noindent\textcolor{Green}{\rule{\textwidth}{1pt}}
	(L5)\hspace{2mm} 1,40 sag an-na ugu 20 sag ki-ta \textit{m}[\textit{i-nam} dirig 1,20 dirig]\\
	(L6)\hspace{2mm} 1,[20] \textit{a-na} 6 \'{e}-gar$_8$ dal-$<$ba$>$-na \'{i}l 8 \textit{he}-[\textit{p\'{e}-ma} 4]
	\\
	(L7)\hspace{2mm} [\textit{ta}]-\textit{ta-a-ar} 1,40 sag an-na k\'{u}-[k\'{u}-\textit{ma} 2,46,40]
	\\
	(L8)\hspace{2mm} [2]0 sag ki-ta k\'{u}-k\'{u}-\textit{ma} 6,[40 \textit{a-na} 2,46,40 \textit{\c{s}\'{i}-ib-ma}]\\
	(L9)\hspace{2mm} 2,53,20 bar 2,53,20-da \'{i}[b-si$_8$-bi \textit{ka-bi-is} 1,12]
	\\
	\noindent\textcolor{Green}{\rule{\textwidth}{1pt}}
	(L10)\hspace{0mm} 1,12 \textit{a-na} \'{e}-gar$_8$ dal-$<$ba$>$-na \textit{a-na} sag [\textit{\v{s}u-ku-un}]
	\\
	(L11)\hspace{0mm} \textit{ta-ta-a-ar} 1,40 sag an-$<$na$>$ ugu [20 sag ki-ta]
	\\
	(L12)\hspace{0mm} 1,20 dirig 1,20 a-r\'{a} 6 8 bar 8 [4 4]
	\\
	(L13)\hspace{0mm} \textit{i-na} 1,12 kud-\textit{ma} 1,8 1,8 \textit{\`{u}} [20 ul-gar-\textit{ma} 1,28]\\
	(L14)\hspace{0mm} bar-\textit{zu} 1,12 \textit{a-na} 6 \'{e}-gar$_8$ dal-$<$ba$>$-[na \'{i}l 7,12]\\
	(L15)\hspace{0mm} 7,12 ki \textit{\v{s}a} \'{e}-gar$_8$ dal-$<$ba$>$-na [$\cdots $ $\cdots $] \\
	(L16)\hspace{0mm} 52,48 bar-\textit{zu} 26,2[4 $\cdots $ $\cdots $ $\cdots $]\\
	(L17)\hspace{0mm} 1,40 \textit{\`{u}} 1,16 [ul-gar-\textit{ma} 2,56 $\cdots $ $\cdots $ $\cdots $]
\end{note1}

\subsection*{Translation}

\underline{Obverse,  Lines 1-15}
\begin{note} 
	\begin{tabbing}
		\hspace{15mm} \= \kill 
		(L1)\> \tabfill{$\cdots $ that of $\cdots $ $\cdots $.}\\
		(L2)\> \tabfill{2,10 of the upper width exceeds 30 of the lower width by 1,40.}\index{width}\\ 
		(L3)\> \tabfill{Make the reciprocal of 3,45 of the length, and (the result is) 0;0,16. Multiply (it) by 1,40, (and the result is) 0;26,40.}\index{length}\index{reciprocal of a number} \\
		(L4)\> \tabfill{2 times 0;26,40 is 0;53,20. Multiply (it) by 4,30,0, (and the result is) 4,0,0.}\\
		(L5)\> \tabfill{You return. Square 2,10 of the upper width, (and the result is) 4,41,40.}\index{width}\\
		(L6)\> \tabfill{Subtract 4,0,0 from 4,41,40, and (the result is) 41,40. (Its square root is 50.)}\index{square root}	
	\end{tabbing}
\end{note} 
\begin{note} 
	\begin{tabbing}
		\hspace{15mm} \= \kill 
		(L7)\> \tabfill{If the middle dividing line is 50, (and) the lower width 30 {\fontfamily{qpl}\selectfont nindan}\index{nindan (length unit)} ($\approx $180m), $\cdots $ $\cdots $.}\index{middle dividing line}\index{width}\\	
		\> \tabfill{\textcolor{gray}{\rule{13.7cm}{1pt}}}\\
		(L8)\> \tabfill{How much 50 of the middle dividing line exceeds 30 of the lower width? It exceeds (30) by 20.}\index{middle dividing line}\index{width}\\ 
		(L9)\> \tabfill{20 times 1,20 is 26,40. 2 times (26,40) is 53,20.}\\ 
		(L10)\> \tabfill{You return. Square 50 of the middle dividing line, and (the result is) 41,40.}\index{middle dividing line} \\
		(L11)\> \tabfill{Its square root is paced off. It is 1,4(?).}\index{square root}	
	\end{tabbing}
\end{note} 
\begin{note} 
	\begin{tabbing}
		\hspace{15mm} \= \kill 
		(L12)\> \tabfill{A trapezoid. 35 is the upper width. 5 is the lower width.}\index{width}\index{trapezoid}\\ 
		(L13)\> \tabfill{Two brothers ought to divide (it) equally.}\\	
		\> \tabfill{\textcolor{gray}{\rule{13.7cm}{1pt}}}\\
		(L14)\> \tabfill{35 times 35 is 20,25. 5 times 5 is 25.}\\
		(L15)\> \tabfill{Add (25) to 20,25, and (the result is) 20,50. Halve (it). 10,25. Its square root is 25.}\index{square root}		
	\end{tabbing}
\end{note} 

\noindent
\underline{Reverse,  Lines 1-17}
\begin{note} 
	\begin{tabbing}
		\hspace{17mm} \= \kill 
		(L1)\> \tabfill{$\cdots $ $\cdots $ $\cdots $ 0;36(?) $\cdots $ $\cdots $.}\\
		(L2)\> \tabfill{1;40 is the upper width. 0;20 is the lower width. 1 is the length. $\cdots $ $\cdots $.}\index{length}\index{width}\\
		(L3)\> \tabfill{In the middle of the length I built a party wall, (whose thickness is) 1 {\fontfamily{qpl}\selectfont k\`{u}\v{s}} 6 {\fontfamily{qpl}\selectfont \v{s}u-si} (= 0;6 {\fontfamily{qpl}\selectfont nindan}\index{nindan (length unit)}).}\index{kuz@k\`{u}\v{s} (length unit)}\index{party wall}\index{length}\\
		(L4)\> \tabfill{Give the remainder (of the trapezoid) to two brothers (equally).}\index{trapezoid}\\
		\> \tabfill{\textcolor{gray}{\rule{13.5cm}{1pt}}}\\
		(L5)\> \tabfill{How much 1;40 of the upper width exceeds 0;20 of the lower width? It exceeds (0;20) by 1;20.}\index{width}\\
		(L6)\> \tabfill{Multiply 1;20 by 0;6 of the party wall, (and the result is) 0;8. Break (it in two), and (the result is) 0;4.}\index{party wall}\\
		(L7)\> \tabfill{You return. Square 1;40 of the upper width, and (the result is) 2;46,40.}\index{width}\\ 
		(L8)\> \tabfill{Square 0;20 of the lower width, and (the result is) 0;6,40. Add (it) to 2;46,40, and (the result is)}\index{width}\\
		(L9)\> \tabfill{2;53,20. With half of 2;53,20, its square root is paced off. It is 1;12.}\index{square root}\\
		\> \tabfill{\textcolor{gray}{\rule{13.5cm}{1pt}}}\\
		(L10)\> \tabfill{Put down 1;12 for the party wall, for the width (of it).}\index{width}\index{party wall}\\
		(L11-12)\> \tabfill{You return. 1;40 of the upper width exceeds 0;20 of the lower width by 1;20. 1;20 times 0;6 is 0;8. Half of 0;8 is 0;4.}\index{width}\\ 
		(L13)\> \tabfill{Subtract 0;4 from 1;12, and (the result is) 1;8. Add 1;8 and 0;20 together, and (the result is) 1;28.}\\ 
		(L14)\> \tabfill{A half. Multiply 1;12 by 0;6 of the party wall, (and the result is) 0;7,12.}\index{party wall}\\ 
		(L15)\> \tabfill{0;7,12 is the area of the party wall. $\cdots $ $\cdots $.}\index{party wall}\index{area of a party wall}\\
		(L16)\> \tabfill{0;52,48. Half (of it) is 0;26,24. $\cdots $ $\cdots $.}\\
		(L17)\> \tabfill{Add 1;40 and 1;16 together, and (the result is) 2;56. $\cdots $ $\cdots $.}
	\end{tabbing}
\end{note}

\subsection*{Mathematical Interpretation}
\begin{figure}[H]
	\centering
	\includegraphics[scale=1]{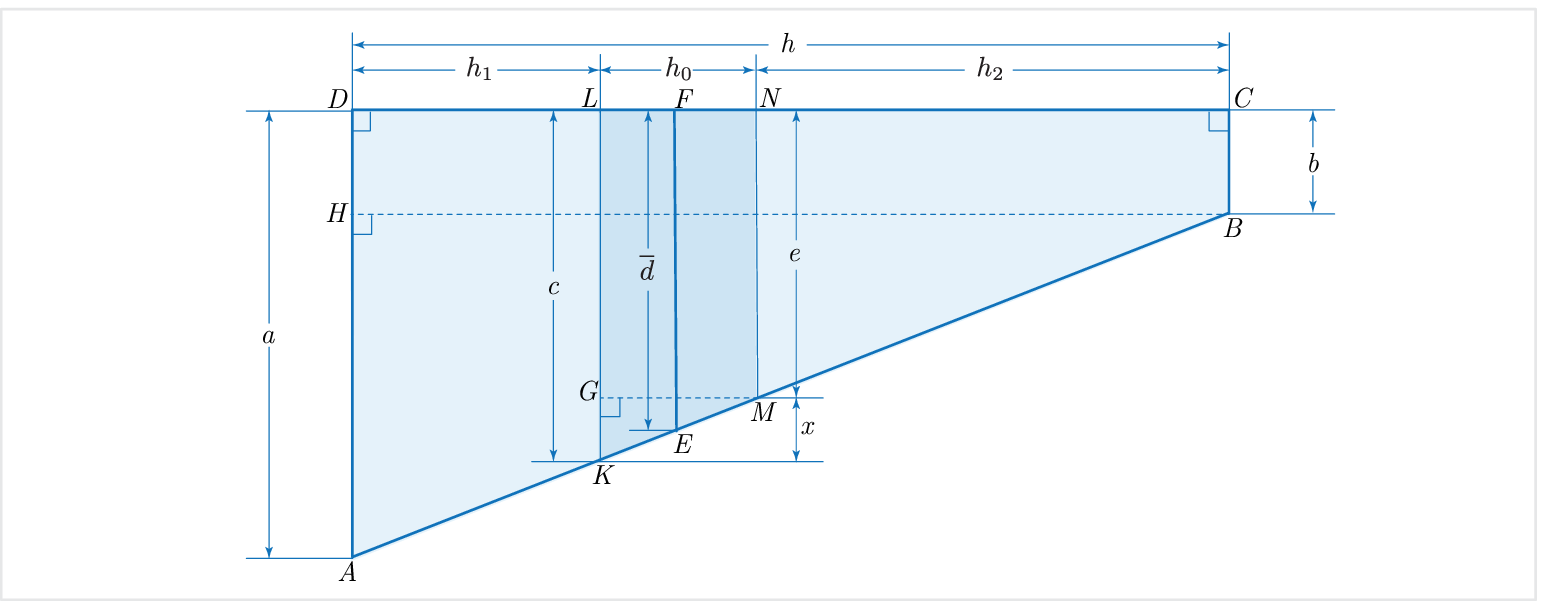}
	\caption{Bisection of a right trapezoid using a party wall}
	\label{Figure11}
\end{figure}

Although   parts of the text regarding the statement of this problem are lost, we can formulate the problem   based on the calculations provided. Consider a right trapezoid  $ABCD$ with bases of lengths $a=1;40$, $b=0;20$ and height  $h=1$ and divide  it into two   parts by   a party wall (the trapezoid  $KMNL$) of width  $h_0=0;6 $ such that   the transversal\index{transversal line} line $EF$  is in the middle of  the left and the right edges $LK$ and $NM$ of the party wall  respectively (see \cref{Figure11}).  The problem now  is  how to determine the values of the edges and the middle line of the party wall  $KMNL$  provided that  the two right trapezoids  $ AKLD$ and $MBCN$ have  equal areas.

Let $\overline{LK}=c $, $ \overline{NM}=e$, $ \overline{FE}=\overline{d}$ and $ \overline{GK}=x$. Clearly, $x=c-e $.  Let $d$ be the bisecting transversal of the trapezoid and $\overline{d} $ be the middle line of the party wall (note that we only have shown $\overline{d}$ in \cref{Figure11}).  The scribe  first computes the difference between the left and the right  edges $LK$ and $NM$ of the party wall and then  finds the length  of the bisecting transversal  line $d$ in lines 5-9. Consider the perpendicular  lines from $M$ onto $LK$ and from $B$ onto $AD$ with intersection points $G$ and $H$ respectively (see \cref{Figure11}). Since two right triangles  $ \triangle MGK$ and $\triangle ABH$ are similar, we have
\[ \frac{\overline{GK}}{\overline{GM}}=\frac{\overline{AH}}{\overline{HB}}  \] 
or equivalently
\[ \frac{x}{h_0}=\frac{a-b}{h}. \] 
This implies that
\begin{equation}\label{equ-m}
	x=\dfrac{h_0}{h}(a-b).
\end{equation}

 It is clear from \cref{Figure11} that the difference between the lengths of the edges of the party wall  and the middle line of the party wall is half the difference between two edges of the party wall. So  this value is  clearly   $ \frac{x}{2}$ and according to lines 5-6, we can use    \cref{equ-m} to  write
 	\begin{align*}
 		\dfrac{x}{2}&=\frac{1}{2}\times\dfrac{h_0}{h}\times(a-b)\\
 		&=(0;30)\times \dfrac{(0;6)}{1} \times(1;40-0;20)\\
 		&=(0;30)\times (0;6)\times(1;20)\\
 		&=(0;30)\times (0;8)\\
 		&=0;4.
 	\end{align*}
 	Hence
 	\begin{equation}\label{equ-n}
 		\dfrac{x}{2}=0;4 
 	\end{equation}
 	or equivalently 
 	\begin{equation}\label{equ-o}
 		x=0;8.
 	\end{equation}
 
 	On the other hand, according to lines 7-9, the value of  bisecting transversal $d$, as usual, is obtained    by  the following calculations: 
 	\begin{align*}
 		d&=\sqrt{\dfrac{a^2+b^2}{2}}\\
 		&=\sqrt{\dfrac{(1;40)^2+(0;20)^2}{2}}\\
 		&=\sqrt{\dfrac{2;46,40+0;6,40}{2}}\\
 		&=\sqrt{\dfrac{2;53,20}{2}}\\
 		&=\sqrt{1;26,40} \\
 		&=1;12,6,39,41,30,\cdots\\
 		&\approx 1;12. 
 	\end{align*}
 	
\begin{figure}[H]
	\centering
	\includegraphics[scale=1]{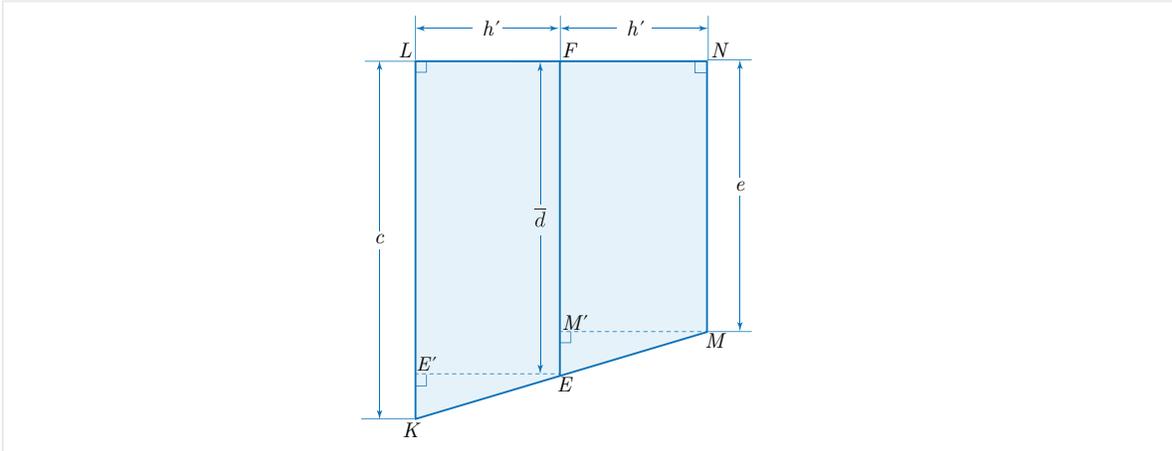}
	\caption{Dividing  a party wall using a transversal line}
	\label{Figure12}
\end{figure}

From now on, the scribe   appears to have  assumed the value  $1;12$ for     $\overline{d}$ and has  utilized it  in the rest of the calculations. So, we set
\begin{equation}\label{equ-p}
	\overline{d}= 1;12.
\end{equation}
Note that  $\overline{d}$ is the mean value of two numbers $c$ and $e$ because by the similarity  of two  right triangles  $\triangle KEE' $ and $ \triangle EMM' $  in \cref{Figure12}, we have 	
 	
\[ \frac{\overline{KE'}}{\overline{EE'}}=\frac{\overline{EM'}}{\overline{MM'}}, \]
since $\overline{KE'}=c-\overline{d} $, $\overline{EM'}=\overline{d}-e$ and $ \overline{EE'}=\overline{MM'}=h'$, the last equality implies that
\[ \frac{c-\overline{d}}{h'}=\frac{\overline{d}-e}{h'}  \]
or
\[  c-\overline{d} = \overline{d}-e. \] 	
This clearly gives the following formula
\begin{equation}\label{equ-q}
	\overline{d}=\frac{c+e}{2}. 
\end{equation}
Maybe this is the reason in line 10 the scribe calls the value $\overline{d}=1;12$ as the ``width'' of the party wall.

The second part of the solution (lines 11-17) deals with finding the areas of two trapezoids   $ADLK $ and $ BCNM$ which is the main goal of the scribe. At the beginning of this part (lines 11-12)  in our opinion the same calculations as in lines 5-6 are perhaps   repeated. Moreover, the calculations in lines 13-16 are  (to modern eyes) somewhat  disorderly. In fact,   lines 14 and 15 calculate  the area of the base of the party wall, i.e., the trapezoid  $KMNL$, and   lines 13 and 16 the area of the trapezoid $ MBCN$.  For the area  of the party wall $KMNL$, say $S_0$, we note that since $\overline{d}$ is the mean value of those of   bases of  the trapezoid  $KMNL$  and its height is $h_0=0;6$, it follows from \cref{equ-p}  and \cref{equ-q}   that
\begin{align*}
	S_{0}&=\frac{h_0(c+e)}{2}\\
	& = h_0\overline{d}\\
	& =(0;6)\times (1;12)\\
	&=0;7,12.
\end{align*}
Thus
\begin{equation}\label{equ-r}
	S_{0}= 0;7,12.
\end{equation}
For  the area  of the trapezoid  $ MBCN$, we need to find  the values of  bases $c,e$ and the height  $h_2$.    It is obvious from \cref{Figure12}  that
\begin{equation}\label{equ-s}
	\begin{cases}
		c=\overline{d}+\frac{x}{2} \\
		e=\overline{d}-\frac{x}{2}.
	\end{cases}
\end{equation}
So, it follows from \cref{equ-o}, \cref{equ-p} and \cref{equ-s} that
\begin{equation}\label{equ-t}
	c=1;12+0;4=1;16 
\end{equation}
and
\begin{equation}\label{equ-u}
	e=1;12-0;4=1;8.
\end{equation}

Although there is no trace of values  for the two heights $h_1$ and $h_2$  in the text, it seems that the scribe was aware that the correct values of  these two heights  are 0;18 and 0;36 respectively. In fact, he is implicitly using these numbers  in his calculations. According to lines 13 and 15,  it follows  from   \cref{equ-u}  that
\begin{align*}
	S_{MBCN}&=\frac{1}{2}\times h_2 \times (e+b)\\
	&=\frac{1}{2}\times  (0;36) \times (1;8+0;20)\\
	&=\frac{1}{2}\times  (0;36) \times (1;28)\\
	&=\frac{1}{2}\times  (0;52,48)  \\
	&=0;26,24. 
\end{align*}
Hence
\begin{equation}\label{equ-SMT26-m}
	S_{MBCN}=0;26,24.
\end{equation} 
In line 17 and the subsequent missing lines, the scribe might have computed the area  of the other trapezoid  $AKLD $ by implicitly using the correct value  $h_1=0;18$. As a matter of fact, it follows from \cref{equ-t}  that
\begin{align*}
	S_{AKLD}&=\frac{1}{2}\times h_1 \times (c+a)\\
	&=\frac{1}{2}\times  (0;18) \times (1;16+1;40)\\
	&=\frac{1}{2}\times  (0;18) \times (2;56)\\
	&=\frac{1}{2}\times  (0;52,48)  \\
	&=0;26,24 
\end{align*}
hence
\begin{equation}\label{equ-v}
	S_{AKLD}=0;26,24.
\end{equation} 
Therefore, the equal share of each brother in land is $ 0;26,24 $.

\section{Significance of \textbf{SMT No.\,26}}
Although one might think that the areas of the two trapezoids   $ MBCN$ and $AKLD$  given   in \cref{equ-u} and \cref{equ-v} respectively are approximate  values because the scribe  has used an approximate value in \cref{equ-q} for $\overline{d} $,   these area values  surprisingly    turn out to be  precise! What is really happening here? Is this a  just a lucky approximation or  did the scribe already know that  these are the correct values?  In fact, from \cref{equ-r}, \cref{equ-u} and \cref{equ-v},   we  suggest  that   the total sum of areas of the three parts is
\begin{align*}
	S_{MBCN}+S_{AKLD}+S_0&=2\times (0;26,24)+ 0;7,12\\
	&=0;52,48+0;7,12\\
	&=1,
\end{align*}
which is exactly the area  of the whole trapezoid  $ABCD$.  This suggests: 
\begin{itemize}
	\item[(1)]  the scribe   successfully accomplished his task  of dividing     land between two brothers equally by means of a party wall;  and  
	\item[(2)]    he  seemed to  already know the exact value  of  $ \overline{d}$, $c$, $e$ and the areas of the two   trapezoids   $ MBCN$ and $AKLD$. 
\end{itemize}

Another issue concerns   the values of  heights $h_1$ and $h_2$ of the  two   trapezoids   $ MBCN$ and $AKLD$. How did the scribe  know   the correct values of the  two heights  $h_1$ and $h_2$?  If we   look  into   this problem, we notice that   $a=5b $, $r=5$ and $n=10$, which  are the data in the sixth row of  \cref{Table1}. In other words, the Susa scribe must have known that these numbers would provide   exact areas for  the two   trapezoids   $ MBCN$ and $AKLD$. But how did he know that? Our guess is that he might have designed the problem by using a similar approach to the one we discussed in section 5. However, since finding   acceptable values for $a,b$ and $n$ might have been an arduous task in the sexagesimal numeral system, he might  have utilized  small values for $n$ and $r$ (for example less than 60) and then performed some calculation to find a correct answer. With some luck, he could have found the correct answers for  $n$ and $r$.

 As    he had to    divide the height  into $n$ equal parts, the value of $n$ must  have been  a regular\footnote{Any number in the form of $2^{p}3^{q}5^{r}$, where $p,q,r$ are nonnegative integers, is called a regular number.} number. By looking at the data in  \cref{Table1},   those  values for     $r$ that provide answers are  2,3,4,5,6,8,9 and the corresponding values for $n$   are 37,17,305,65,10,25,35,20 among which only 20, 25 and 20 are regular numbers. Although, the acceptable pairs of $(r,n)$ are $(5,10), (6,25)$ and $(9,20)$,   the first acceptable pair    $ (5,10) $ has the minimum value for $n$. One might reasonably anticipate     this acceptable pair of values for $r,n$  to be the obvious choice and  surprisingly enough, these numbers  are the very values that the Susa scribe has used in this problem! Using either the approach we discussed or through trial and error, the Susa scribe must have known these numbers in order to design this problem. Without using the correct choices for $r$ and $n$, the calculations would be very complicated and the answers   approximate not exact.

It is   reasonable   to    posit  that the Susa scribe  of this text   who designed such a beautiful and elegant  problem  employed an innovative approach by which  means    all the exact values for the required parts in   question were obtained.

\section{Conclusion}
Like the Babylonians, the  Susa scribes  dealt with transversal bisectors of trapezoids and solved   problems using the Babylonian formula for the transversal. However, the Susa scribe of \textbf{SMT No.\,26} dared to consider a more general case and instead of transversal lines he   dealt with transversal strips. Although one may think treating such a problem is similar to a transversal line case, the reality is different and this problem requires a higher level of  mathematical skill and experience. It seems that the Susa scribe was cognizant of   the numbers that would lead   to the solutions. These numbers   could have been determined  only by finding  the natural solutions of a quadratic equation whose coefficients are with respect to the bases and the number of strips. In fact, the only way to find such  solutions is to  check  a huge number of different choices for both the bases and    the number of strips, which requires performing complicated calculations in the sexagesimal numeral system.

The  design of  mathematical problems requires  deep knowledge and a great deal of experience  on the part of a teacher. A beautiful and challenging problem  can  engage the   imagination  of students and  encourage  their creativity. The mathematical interpretation of the fourth problem in \textbf{SMT No.\,26} that we have identified in this article  displays  just such  a creative approach to  both the design and the means of  solving the problem.   Among the many and varied   characteristics attributable  to  the Susa scribes apparent from the \textbf{SMT},  we consider  this  characteristic  to be   by far   the most striking one.

\end{document}